# A Two-Layer Distributed Control Method for Islanded Networked Microgrid Systems

Xiaoyu Wu, *Student Member, IEEE*, Yin Xu, *Member, IEEE*, Xiangyu Wu, *Member, IEEE*, Jinghan He, *Senior Member, IEEE*, Josep M. Guerrero, *Fellow, IEEE*, Chen-Ching Liu, *Fellow, IEEE*, Kevin P. Schneider, *Senior Member, IEEE* and Dan T. Ton

*Abstract*—**This paper presents a two-layer, four-level distributed control method for networked microgrid (NMG) systems, taking into account the proprietary nature of microgrid (MG) owners. The proposed control architecture consists of a MG-control layer and a NMG-control layer. In the MG layer, the primary and distributed secondary control realize accurate power sharing among distributed generators (DGs) and the frequency/voltage reference following within each MG. In the NMG layer, the tertiary control enables regulation of the power flowing through the point of common coupling (PCC) of each MG in a decentralized manner. Furthermore, the distributed quaternary control restores system frequency and critical bus voltage to their nominal values and ensures accurate power sharing among MGs. A small-signal dynamic model is developed to evaluate dynamic performance of NMG systems with the proposed control method. Time-domain simulations as well as experiments on NMG test systems are performed to validate the effectiveness of the proposed method.**

## I. INTRODUCTION

R esiliency against major disasters, such as major hurricanes or earthquakes, is considered by the U.S. Department of Energy (DOE) as a most essential characteristic of the future smart distribution systems [1]. Concerning the enhancement of system resiliency, interconnecting microgrids (MGs) to form a networked microgrid (NMG) system after a major outage has been proved an effective option [2].

Three types of NMG systems (or multi-microgrid systems) are reported in the literature: low voltage (LV) MGs interconnected through LV tie lines [3], medium voltage (MV) MGs networked via MV feeders [4], and LV MGs interconnected through a MV feeder and distribution transformers [5], [6].

This work was supported by Natural Science Foundation of China (51807004, 51807005) and the U.S. Department of Energy (DOE). (*Corresponding author: Xiangyu Wu*)

Xiaoyu Wu, Yin Xu, Xiangyu Wu, Jinghan He are with the School of Electrical Engineering, Beijing Jiaotong University, Beijing 100044, China.

Josep M. Guerrero is with the Department of Energy Technology, Aalborg University, 9220 Aalborg East, Denmark.

Chen-Ching Liu is with the Department of Electrical and Computer Engineering, Virginia Polytechnic Institute and State University (VT), Blacksburg, VA 24061, USA.

Kevin P. Schneider is with the Pacific Northwest National Laboratory (PNNL) located at the Battelle Seattle Research Center in Seattle, Washington.

Dan T. Ton is with the U.S. Department of Energy (DOE) Office of Electricity Delivery and Energy Reliability (OE), Washington, DC 20585, USA.

Disregarding the way in which MGs are networked, NMG systems have some features in common. Firstly, MGs within a NMG system may belong to different entities, and limited information may be shared with others due to the proprietary nature. Secondly, multiple control objectives need to be realized by effective coordinations among DGs and MGs.

### A. Literature Review

#### 1) Control architectures of NMG systems

A proper control architecture for NMG systems considering the above features is needed. The classical three-level control architecture [7] is widely applied to a single MG [8]-[12]. The existing control philosophies for NMG systems are usually based on this architecture and fall into two categories, i.e., one-layer architecture and two-layer architecture. The one-layer architecture ignores the boundary of each MG, and considers the NMG system as a large MG. Thus, the three-level control architecture of MGs can be modified to fit the NMG system [13]-[15]. The two-layer architecture adds an extra NMG-control layer and considers the three-level control of MGs as the MG-control layer to form a two-layer architecture [16]-[29]. The two-layer architecture is inspired by the multi-layer and multi-area control concept in the bulk power system [30]. The NMG-layer regards each MG as a control entity and thus avoids a direct control of each DG unit in the one-layer architecture. Therefore, the two-layer architecture is preferred in this study.

However, the existing two-layer architectures have the following limitations: (i) they cannot realize load sharing among MGs automatically and enable the plug-and-play capability of each MG, (ii) the NMG-control layer requires too much information of DGs within MGs, which may be inaccessible, e.g. DG capacity information and load consumption data.

#### 2) Control methods of NMG systems

Centralized or distributed methods can both be applied to the control of NMG systems. The centralized control methods for NMG systems are reported in [17]-[24]. In [17], the NMG-layer is responsible for calculating power and voltage reference values and then send them to the MG-layer to realize power sharing and voltage control objectives. In [18], the frequency control issue is dealt with in the NMG-layer by coordinating each MG. In [19-24], an interface converter is assumed to be deployed with each MG to realize power sharing among MGs, thus the NMG-layer is actually responsible for the control of the interfaced converter.

Compared with centralized methods, the distributed control methods have some advantages, e.g., robust to communication failures and good re-configurability due to their neighboring communication features. Therefore, the distributed control methods attract much attention in recent years.

Among the existing publications with distributed control methods [13], [14], [25]-[27], both one-layer architecture [13], [14] and two-layer architecture are adopted [25]-[27] . In [13], [14], a distributed communication network including all DG units of the NMG system exists under the one-layer architecture. Thus, there will be a large and complex communication network, which may result in a slow convergence speed of the distributed control algorithm. In [25], [26], the NMG-layer is distributed, while the MG-layer is still centralized. Thus, the DG units cannot flexibly plug in or out due to the centralized MG-layer. In [27], both MG and NMG layer employ distributed control method. However, the MG output power sharing and the critical bus voltage control objectives cannot be realized.

In sum, the above methods cannot simultaneously realize the control objectives of both NMG layer (frequency/voltage regulation and power sharing among MGs) and MG layer (load sharing among DGs), especially under the two-layer distributed architecture.

*3) Stability modeling and analysis of NMG systems*

Compared with a single MG, the dynamic interactions among MGs and among multiple control layers in a NMG system may introduce new low-damping oscillation modes which may even destabilize the system. Therefore, a small-signal dynamic model of the NMG system and its corresponding stability analysis are of significant importance. However, only several publications about the small-signal stability issues of NMG system are reported [23], [24], [31]. In [23], [24], a small-signal dynamic modeling method for the NMG system is proposed, in which each MG is simplified as a DG unit without considering its internal dynamics. This simplification will inevitably lead to analysis errors. In [31], a detailed small-signal dynamic model of a PV-based NMG is proposed and the analysis results indicate that the coupling among MGs will weaken the system stability margin. However, only decentralized primary control is employed with each DG in [31], which means the impacts of distributed control methods and other control layers are not studied in [31].

Based on the above analysis, to the best of the authors' knowledge, a detailed small-signal dynamic model and corresponding stability analysis of the NMG system considering distributed control methods and multiple control layers have not been reported before.

### B. Contribution

Compared with the state-of-the-art, the major contributions of this paper are three-folds:

1) A two-layer, four-level distributed control architecture is proposed. In this architecture, each MG is represented by an agent and only the total spare power capacity information is provided to the NMG layer. Thus, proprietary information of MG entities is well protected. Besides, an interface level is designed in the NMG layer, which can realize load sharing among MGs automatically as well as enable the plug-and-play capability of each MG.

2) In the NMG-control layer, the control method for the interface level is proposed and a distributed control strategy based on it is developed. In the MG-control layer, the classical MG distributed control is adjusted and accommodated with the NMG-layer. The proposed control strategy is capable of simultaneously i) regulating system frequency and critical bus voltage to desired values, and ii) achieving accurate active and reactive power sharing among MGs as well as among DGs within each MG.

3) A unified small-signal dynamic model of the NMG system with the proposed control strategy is constructed. The model is detailed enough, which means the dynamics of every line, load, especially the two-layer distributed controllers are included. Moreover, the stability analysis based on the proposed model reveals the newly introduced low-damping oscillation modes and their impact factors. Then, general guidelines are provided for controller parameters tuning based on the stability assessment results.

### C. Paper Organization

The remainder of this paper is organized as follows. The proposed control architecture is introduced in Section II. In Section III, the design of controllers and corresponding coordination principles are presented. Section IV develops a unified small-signal dynamic model of a test NMG system. Stability analysis and numerical simulation results are discussed in Section V. Section VI provides the experimental results. Conclusions are summarized in Section VII.

## II. THE HIERARCHICAL CONTROL ARCHITECTURE

In this section, the control objectives and overall control architecture for NMG systems are described.

### A. Control Objectives

This paper presents a hierarchical control architecture to perform frequency and voltage regulation and power sharing control of the NMG system. The control objectives include:

(i) To maintain the system frequency $f_{sys}$ at its rated value $f_{sys}^*$.

(ii) To restore the critical bus voltage $V_c$ to the desired value $V_c^*$. Note that only one critical bus is assumed, and it can be selected according to the operation requirement.

(iii) To share active and reactive power among MGs in proportion to the reference values, i.e.,

$$\frac{1}{P_{PCC1}^*}P_{PCC1} = \frac{1}{P_{PCC2}^*}P_{PCC2} = \cdots = \frac{1}{P_{PCCm}^*}P_{PCCm} \quad (1)$$

$$\frac{1}{Q_{PCC1}^*}Q_{PCC1} = \frac{1}{Q_{PCC2}^*}Q_{PCC2} = \cdots = \frac{1}{Q_{PCCm}^*}Q_{PCCm} \quad (2)$$

where $P_{PCCk}^*$, $Q_{PCCk}^*$, $P_{PCCk}$, $Q_{PCCk}$ are active and reactive power reference for point of common coupling (PCC) of MG$_k$, output active and reactive power of PCC$_k$, respectively, with $k \in \mathcal{M}$, $\mathcal{M} = \{1, 2, ..., m\}$. In this paper, $P_{PCCk}^*$, $Q_{PCCk}^*$ are equal to the spare active and reactive capacity of MG$_k$, denoted by $P_{sMGk}$, $Q_{sMGk}$, respectively. Note that $P_{PCCk}^*$, $Q_{PCCk}^*$ can also be determined according to optimal power flow results, which makes the control objectives flexible. Equations (1) and (2) are denoted as objective (iii)-(1) and objective (iii)-(2),

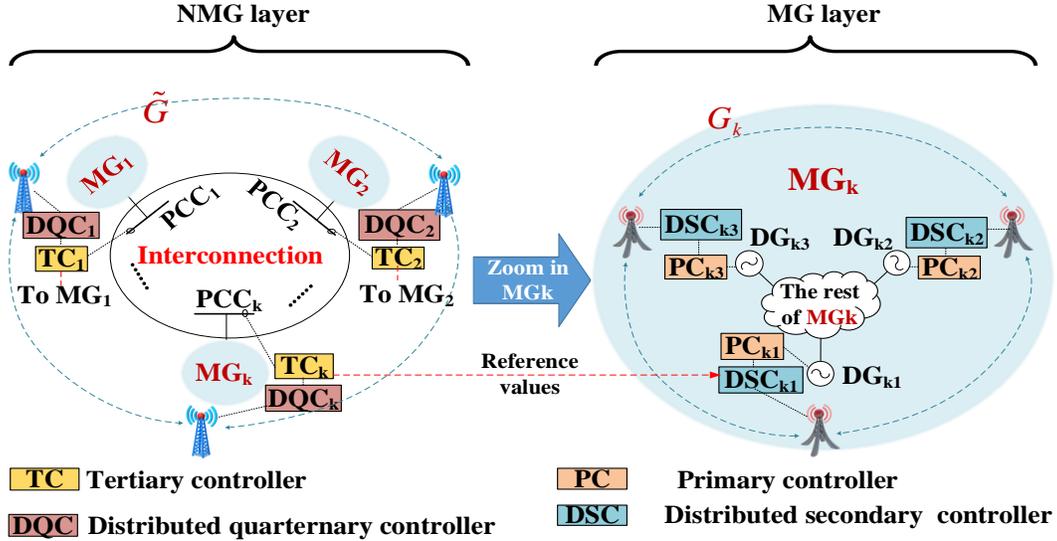

Fig. 1  A schematic diagram of a NMG system with the proposed strategy

respectively.

(iv) Within each MG, active and reactive power can be shared among DGs based on their power capacities, i.e.,

$$\frac{1}{P_{\max k1}}P_{k1} = \frac{1}{P_{\max k2}}P_{k2} = \cdots = \frac{1}{P_{\max kn_k}}P_{kn_k} \qquad (3)$$

$$\frac{1}{Q_{\max k1}}Q_{k1} = \frac{1}{Q_{\max k1}}Q_{k2} = \cdots = \frac{1}{Q_{\max kn_k}}Q_{kn_k} \qquad (4)$$

where $P_{\max ki}$, $Q_{\max ki}$, $P_{ki}$, $Q_{ki}$ are active and reactive power capacities, active and reactive power outputs of $DG_i$ in $MG_k$, respectively, with $i \in \mathcal{N}_k$, $\mathcal{N}_k = \{1,2,\dots,n_k\}$. Equations (3) and (4) are denoted as objective (iv)-(1) and objective (iv)-(2), respectively.

In the authors' previous work [6], a NMG power flow model considering the above objectives is proposed and only one solution exists, which demonstrates that the objectives (i)~(iv) can be met simultaneously.

### B. The Proposed Two-Layer Control Architecture

In order to realize the aforementioned control objectives, a control architecture, shown in Fig. 1, which consists of two layers with four levels is proposed. From the bottom to the top, they are primary and distributed secondary levels in the MG-control layer, tertiary and distributed quaternary levels in the NMG-control layer, respectively. The primary level is responsible for DG's local power, voltage, and current control. Then, the distributed secondary level is introduced to maintain the voltage and frequency of each MG at the reference values received from the tertiary level. The tertiary level manages the power flow of PCC by sending commands to the secondary level. At the top of this architecture, the distributed quaternary level supervises the entire NMG system and control the critical point voltage and system frequency at desired values.

*1) MG-control layer:* This layer aims at meeting the frequency, voltage and power sharing control objectives of MG layer as well as supporting the NMG layer control.

a) The *primary controller (PC) level* is responsible for

regulating the output power of DGs via the droop method [11].

b) The *distributed secondary controller (DSC) level* is responsible for regulating the MG frequency and PCC voltage according to the reference values received from the tertiary controller level. The control actions are taken by sending compensation signals to the primary controllers. A distributed communication network $G$ is set for each MG, and DGs inside the MG are assumed to be communication nodes.

*2) NMG-control layer:* This layer handles the control of system frequency and critical bus voltage as well as power sharing among MGs.

c) The *tertiary controller (TC) level* is responsible for sharing power among MGs by controlling the PCC power flow according to the droop characteristics. Through adjusting the frequency and PCC voltage reference values which are sent to the distributed secondary level, the PCC power flow is controlled. Note that only the total spare capacity information of MGs is needed in this level.

d) The *distributed quaternary controller (DQC) level* regulates the system frequency and critical bus voltage to their desired values by coordinating MGs in a distributed manner. Each MG is modeled as an agent to form the distributed communication network $\tilde{G}$ in this level.

Note that the NMG layer communication network $\tilde{G}$ needs a connection with every $G$ in the MG-control layer, as shown by the red dotted link in Fig. 1.

### III. CONTROLLER DESIGN

Based on the proposed architecture, the corresponding controllers and coordination strategy are described in this section. Fig. 2 shows the control block diagram of a NMG system which consists of $m$ LV MGs interconnected through a MV feeder. Note that the proposed control strategy is also suitable for other

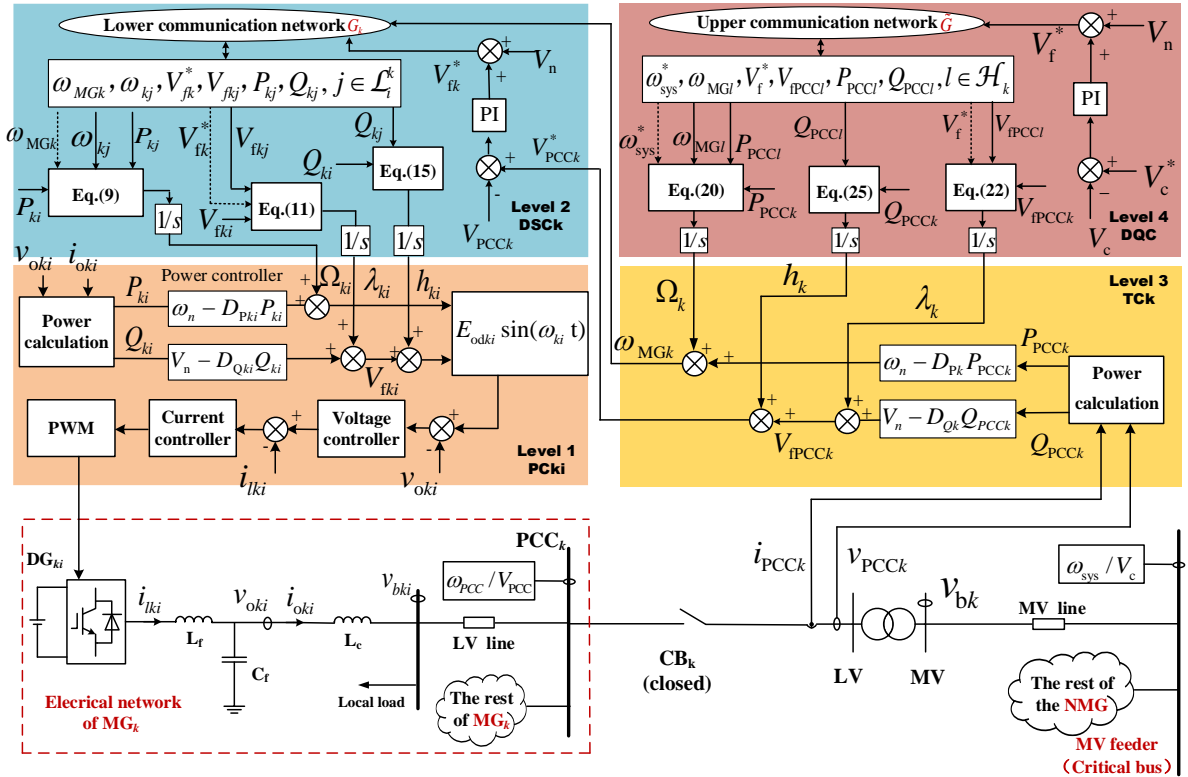

Fig. 2  A block diagram of the proposed two-layer distributed control method

types of NMG system introduced in Section I.

### A. PC Level

In this level, the droop-based control is adopted, which consists of the power controller, inner voltage controller and current controller, as shown in Fig. 2. The power controller allows DGs to share active and reactive power demand based on their power capacities by setting droop coefficients, i.e.,

$$\omega_{ki} = \omega_n - D_{Pki}P_{ki} \tag{5}$$

$$V_{fki} = V_n - D_{Qki}Q_{ki} \tag{6}$$

where $\omega_{ki}$ is the angular frequency of $DG_i$ in $MG_k$, $\omega_n$ is the rated angular frequency, $V_n$ is the rated voltage of the LV network, and $V_{fki}$ is the inverter AC-side voltage reference provided to the inner voltage controller. $D_{Pki}$ and $D_{Qki}$ are the active and reactive power droop coefficients, given by

$$D_{Pki} = \frac{\omega_{max} - \omega_{min}}{P_{maxki}}, \quad D_{Qki} = \frac{V_{max} - V_{min}}{Q_{maxki}} \tag{7}$$

where $\omega_{max}$ and $\omega_{min}$ are the upper and lower limits of the angular frequency, respectively. $V_{max}$ and $V_{min}$ are the upper and lower limits of the DG output voltage, respectively.

### B. DSC Level

The DSC level is responsible for realizing power sharing objectives of DGs within each MG as well as tracking the voltage and frequency reference values from the tertiary level. The method introduced in this level is based on the continuous-time average consensus algorithm [32]. The communication network in this level contains $m$ digraphs, $G_1, G_2 \ldots G_m$, corresponding to $m$ MGs, respectively. For $MG_k$, each DG is considered as a node in $G_k$. For DG node $i$, the set of neighbors is denoted as $\mathcal{L}_i^k$. Each node requires its own information and that of its neighbor $j$ ($j \in \mathcal{L}_i^k$) on the digraph to update its states. The associated adjacency matrix is $\mathbf{A}^k = [a_{ij}^k]$, and $a_{ij}^k$ is the weight of edge between node $i$ and $j$.

1) *Distributed secondary frequency control:* With this control, objectives (iv)-(3) can be achieved. Besides, the frequency reference value $\omega_{MGk}$ from $TC_k$ can be tracked. The reference value $\omega_{MGk}$ for different MGs can be different during transient events to adjust the PCC power flow, but will converge to the system rated angular frequency $\omega_{sys}^*$ when the system reaches a steady state. The controller design is a combination of the regulator synchronization problem [32] and the tracking synchronization problem [33], given by

$$\omega_{ki} = \omega_n - D_{Pki}P_{ki} + \Omega_{ki} \tag{8}$$

$$\frac{d\Omega_{ki}}{dt} = -c_{\omega ki}\left[\sum_{j\in\mathcal{L}_i^k} a_{ij}^k(\omega_{ki}-\omega_{kj}) + g_i^k(\omega_{ki}-\omega_{MGk})\right]$$
$$- c_{pki}\sum_{j\in\mathcal{L}_i^k} a_{ij}^k(D_{Pki}P_{ki}-D_{Pkj}P_{kj}) \tag{9}$$

where $c_{\omega ki}$ and $c_{pki}$ are the positive control gains, the pinning gain $g_i^k \geq 0$ is the weight of edge connected to the reference. It is non-zero only for a few nodes (at least one node). Equation (8) is transformed from (5) with an additional DSC variable $\Omega_{ki}$.

2) *Distributed secondary PCC voltage control:* This controller is responsible for controlling each MG's PCC voltage to the reference $V_{PCCk}^*$ from $TC_k$. The correction term $\lambda_{ki}$ is added in the reactive power droop control (6), i.e.,

$$V_{fki} = V_n - D_{Qki}Q_{ki} + \lambda_{ki} \tag{10}$$

$$\frac{d\lambda_{ki}}{dt} = -c_{vki}[\sum_{j\in\mathcal{L}_i^k} a_{ij}^k(V_{fki}-V_{fkj}) + g_i^k(V_{fki}-V_{fk}^*)] \tag{11}$$

where $c_{vki}$ is a positive control gain. $\sum_{j\in\mathcal{L}_i^k} a_{ij}^k(V_{fki}-V_{fkj}) + g_i^k(V_{fki}-V_{fk}^*)$ is the local neighbor tracking error of $V_{fki}$, which enables voltage regulation. $V_{\mathrm{f}k}^*$ is generated through a PI controller such that $V_{\mathrm{PCC}k}$ recovers to its reference $V_{\mathrm{PCC}k}^*$ which is received from $TC_k$, i.e.,

$$V_{\mathrm{f}k}^* = V_{\mathrm{n}} + k_{\mathrm{p}k}(V_{\mathrm{PCC}k}^* - V_{\mathrm{PCC}k}) + k_{\mathrm{i}k}\int(V_{\mathrm{PCC}k}^* - V_{\mathrm{PCC}k})\,dt \quad (12)$$

where $k_{\mathrm{p}k}$ and $k_{\mathrm{i}k}$ are the gains of the PI controller.

*3) Distributed secondary reactive power control:* This controller deals with the inaccuracy of reactive power sharing problem due to the unbalanced line impedance [10]. Thus, the voltage correction term $h_{ki}$ is added to the right-hand side of (10) to realize accurate reactive power sharing among DGs within $MG_k$, namely objective (iv)- (4), by regulating the voltage reference, i.e.,

$$E_{\mathrm{od}ki} = \overbrace{V_{\mathrm{n}} - D_{Qki}Q_{ki} + \lambda_{ki}}^{V_{fki}} + h_{ki} \quad (13)$$

$$E_{\mathrm{oq}ki} = 0 \quad (14)$$

where the voltage reference, i.e., $V_{\mathrm{n}} - D_{Qki}Q_{ki} + \lambda_{ki} + h_{ki}$, is aligned to the d-axis of the $DG_{ki}$ local dq-frame, and the q-axis reference, i.e., $E_{\mathrm{oq}ki}$ is set to zero. Then, $E_{\mathrm{od}ki}$ and $E_{\mathrm{oq}ki}$ are provided to the inner voltage controller. $h_{ki}$ is selected such that $D_{Qki}Q_{ki}$ of each DG in $MG_k$ converges to a common value, which is a regulator synchronization problem [32] given by

$$\frac{dh_{ki}}{dt} = -c_{\mathrm{q}ki}\sum_{j\in\mathcal{L}_i^k}a_{ij}^k(D_{Qki}Q_{ki} - D_{Qkj}Q_{kj}) \quad (15)$$

where $c_{\mathrm{q}ki}$ is a positive control gain, and $\sum_{j\in\mathcal{L}_i^k}a_{ki}(D_{Qki}Q_{ki} - D_{Qkj}Q_{kj})$ is the local neighbor tracking error which enables accurate reactive power sharing.

### C. TC Level

This level is an interface level which can realize load sharing among MGs automatically as well as enable the plug-and-play capability of each MG. Besides, the proprietary information of each MG can be well protected by introducing this level. The droop control is modified to control the output power through PCC of MGs, i.e.,

$$\omega_{\mathrm{MG}k} = \omega_n - D_{\mathrm{P}k}P_{\mathrm{PCC}k} \quad (16)$$

$$V_{\mathrm{fPCC}k} = V_n - D_{Qk}Q_{\mathrm{PCC}k} \quad (17)$$

where $D_{\mathrm{P}k}$ and $D_{Qk}$ are the active and reactive droop coefficients of $MG_k$, respectively, determined by

$$D_{\mathrm{P}k} = \frac{\omega_{\max} - \omega_{\min}}{P_{\mathrm{SMG}k}}, \quad D_{Qk} = \frac{V_{\max} - V_{\min}}{Q_{\mathrm{SMG}k}} \quad (18)$$

In this level, each MG is considered as a droop-controlled node and only the spare power capacity of each MG is needed.

### D. DQC Level

The DQC level is responsible for regulating system frequency and critical bus voltage to desired values, as well as realizing accurate power sharing among MGs. The communication network in this level is denoted as $\tilde{G}$ with associated adjacency matrix $\mathbf{A} = [a_{kl}]$. Each MG is a node in $\tilde{G}$. For MG node $k$, the set of neighbors of node $k$ is denoted as $\mathcal{H}_k$. Three controllers are included in this level.

*1) Distributed quaternary frequency control:* With this control, objective (i) and objective (iii)-(1) can be achieved. The controller design is as follows.

$$\omega_{\mathrm{MG}k} = \omega_{\mathrm{n}} - D_{\mathrm{P}k}P_{\mathrm{PCC}k} + \Omega_k \quad (19)$$

$$\frac{d\Omega_k}{dt} = -c_{\omega k}\left[\sum_{l\in\mathcal{H}_k}a_k(\omega_{\mathrm{MG}k} - \omega_{\mathrm{MG}l}) + g_k(\omega_{\mathrm{MG}k} - \omega_{\mathrm{sys}}^*)\right] - c_{\mathrm{p}k}\sum_{l\in\mathcal{H}_k}a_k(D_{\mathrm{P}k}P_{\mathrm{PCC}k} - D_{\mathrm{P}l}P_{\mathrm{PCC}l}) \quad (20)$$

where $c_{\omega k}$ and $c_{\mathrm{p}k}$ are positive control gains. Equation (19) is (16) with an additional quaternary control variable $\Omega_k$. $\omega_{\mathrm{MG}k}$ is generated as the frequency reference of (9).

*2) Distributed quaternary critical bus voltage control:* This controller is responsible for achieving objective (ii). The term $\lambda_k$ is added in MG's reactive power droop controller (17), i.e.,

$$V_{\mathrm{fPCC}k} = V_{\mathrm{n}} - D_{Qk}Q_{\mathrm{PCC}k} + \lambda_k \quad (21)$$

$$\frac{d\lambda_k}{dt} = -c_{vk}\left[\sum_{l\in\mathcal{H}_k}a_{kl}(V_{\mathrm{fPCC}k} - V_{\mathrm{fPCC}l}) + g_k(V_{\mathrm{fPCC}k} - V_{\mathrm{f}}^*)\right] \quad (22)$$

where $c_{vk}$ is a positive control gain, $V_{\mathrm{f}}^*$ is generated through a PI controller such that $V_c$ recovers to its reference $V_c^*$, i.e.,

$$V_{\mathrm{f}}^* = V_{\mathrm{n}} + k_{\mathrm{p}}(V_c^* - V_c) + k_{\mathrm{i}}\int(V_c^* - V_c)\,dt \quad (23)$$

where $k_{\mathrm{p}}$ and $k_{\mathrm{i}}$ are the gains of the PI controller.

*3) Distributed quaternary reactive power control:* The inaccurate reactive power sharing among MGs is managed by this controller. The voltage correction term $h_k$ is added to the right side of (21) to achieve accurate reactive power sharing by regulating the voltage reference, i.e.,

$$V_{\mathrm{PCC}k}^* = \overbrace{V_{\mathrm{n}} - D_{Qk}Q_{\mathrm{PCC}k} + \lambda_k}^{V_{\mathrm{fPCC}k}} + h_k \quad (24)$$

where the voltage reference $V_{\mathrm{PCC}k}^*$ is generated as the reference of (12). $h_k$ is selected such that $n_k Q_{\mathrm{PCC}k}$ of each MG converges to a common value given by

$$\frac{dh_k}{dt} = -c_{\mathrm{q}k}\sum_{l\in\mathcal{L}_k}a_{kl}(D_{Qk}Q_{\mathrm{PCC}k} - D_{Ql}Q_{\mathrm{PCC}l}) \quad (25)$$

where $c_{\mathrm{q}k}$ is a positive control gain. Thus, objective (iii)-(2) can be realized.

### E. Flow chart and implementation steps of proposed method

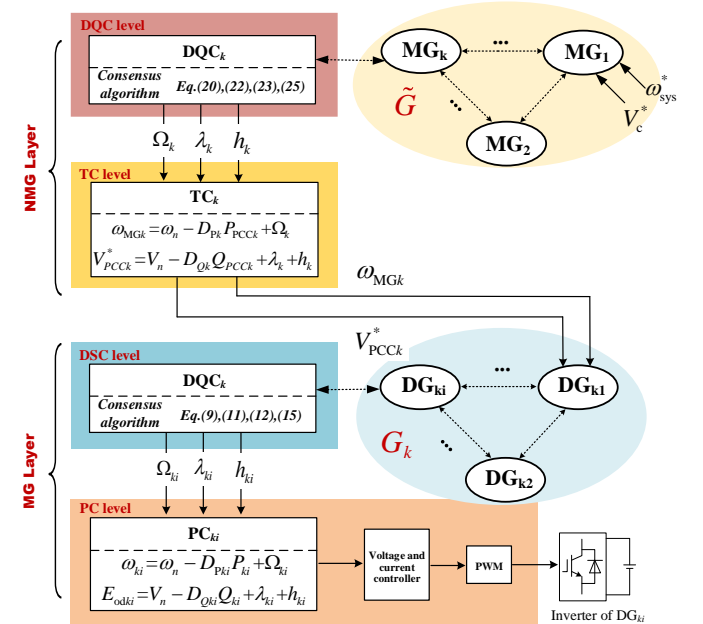

Fig. 3 Flow chart of the proposed control method

In order to illustrate the execution process of the proposed two-layer distributed control method. A flow chart diagram is shown in Fig. 3. Besides, the corresponding implementation steps are given as follows.

*Step 1 (DQC level):* The rated system frequency $\omega_{sys}^*$ and desired critical bus voltage $V_c^*$ are transmitted to agent MG$_l$ in the upper communication network $\bar{G}$. Through the distributed consensus algorithm (20), (22), (23) and (25), the control variables $\Omega_k$, $\lambda_k$ and $h_k$ are obtained and sent to TC level;

*Step 2 (TC level):* The TC is a droop-based controller which is used to adjust PCC power flow of each MG. $\Omega_k$ is applied to shift the frequency-active power droop curve to realize objective (i) and objective (iii)-(1). $\lambda_k$ and $h_k$ are applied to shift the voltage-reactive power droop curve to realize objective (ii) and objective (iii)-(2). The output variables $\omega_{MGk}$ and $V_{PCCk}^*$ are sent to DSC level;

*Step 3 (DSC level):* $\omega_{MGk}$ and $V_{PCCk}^*$ are received by DG$_{kl}$ in the lower communication network $G_k$. Through the distributed consensus algorithm (9), (11), (12) and (15), the control variables $\Omega_{ki}$, $\lambda_{ki}$ and $h_{ki}$ are obtained and sent to PC level;

*Step 4 (PC level):* PC is droop-based for each DG unit. $\Omega_{ki}$ regulates angular frequency of DG$_{ki}$ to $\omega_{MGk}$ and realizes objective (iv)-(1) through shifting the frequency-active power droop curve of PC. $\lambda_{ki}$ and $h_{ki}$ regulate the output voltage of DG$_{ki}$ to $V_{fk}^*$ and realize objective (iv)-(2). The output frequency reference $\omega_{ki}$ and voltage reference $E_{odki}$ are sent to the voltage and current controller, and the switching signals through PWM module are finally generated to control the inverter of DG$_{ki}$.

## IV. SMALL-SIGNAL DYNAMIC MODEL OF NMG SYSTEMS

In order to reveal the dynamic interaction mechanism and guide the controller parameters selection, a small-signal dynamic model of the NMG system in Fig. 2 is developed to evaluate the dynamic performance of the proposed method.

### A. MG Layer Modeling

The MG layer represents the dynamics of the PC and DSC controllers as well as the lines and loads within MGs. Note that each MG is modeled separately and they will be combined in Section IV-C.

*1) DG unit model:* In this paper, the local dq-frame of DG$_{11}$, namely DG$_1$ in MG$_1$, is selected as the common DQ-frame of the entire system. The symbol $\omega_g$ denotes the rotating frequency of DQ-frame and $\omega_g = \omega_{11}$. $\delta_{ki}$ is the angle between the local dq-frame of DG$_{ki}$ and the common DQ-frame, then

$$\dot{\delta}_{ki} = \omega_{ki} - \omega_g \qquad (26)$$

This paper focuses on dynamics of the power controller. Therefore, the relatively fast dynamics of voltage and current controllers can be neglected by assuming

$$v_{odki} = E_{odki}, \ v_{oqki} = E_{oqki} \qquad (27)$$

where $v_{odki}$ and $v_{oqki}$ are the d-axis and q-axis component of DG output voltage $v_{oki}$ as shown in Fig. 2.

By linearizing (8)-(11), (13)-(15) and (26) around an operating point, the model of DG$_{ki}$ can be derived as

$$[\Delta \dot{X}_{DGki}] = A_{DGki}[\Delta X_{DGki}] + B_{DGki}[\Delta v_{bDQki}] + C_{DGki}\Delta \omega_g + \sum_{j \in \mathcal{L}_i^k} F_{DGki}[\Delta X_{DGkj}] + H_{DGki}\Delta V_{fk}^* \qquad (28)$$

where $\Delta v_{bDQki}$ is the deviation of $v_{bki}$ (bus voltage as shown in Fig. 2) in the common DQ-frame, $A_{DGki}$, $B_{DGki}$, $C_{DGki}$, $F_{DGki}$ and $H_{DGki}$ are parameter matrices. Note that $F_{DGki}$ reflects the correlation between DG$_{ki}$ and its neighbors DG$_{kj}$, $j \in \mathcal{L}_i^k$. The state variables of each DG unit are

$$[\Delta X_{DGki}] = [\Delta \delta_{ki}, \Delta P_{ki}, \Delta Q_{ki}, \Delta \Omega_{ki}, \Delta \lambda_{ki}, \Delta h_{ki}, \Delta i_{odki}, \Delta i_{oqki}]^T \qquad (29)$$

*2) PCC voltage controller model:* Introduce $\psi_k$ as the state variable to describe the dynamics of (12), i.e.,

$$\dot{\psi}_k = V_{PCCk}^* - V_{PCCk} \qquad (30)$$

where $V_{PCCk} = \sqrt{V_{PCCDk}^2 + V_{PCCQk}^2}$. Then, the small-signal dynamic model of the PCC voltage controller is obtained through linearizing (30) and (12), i.e.,

$$\Delta \dot{\psi}_k = -A_{PCCk}[\Delta V_{PCCDQk}] + \Delta V_{PCCk}^* \qquad (31)$$

$$\Delta V_{fk}^* = -k_{pk}A_{PCCk}[\Delta V_{PCCDQk}] + k_{pk}\Delta V_{PCCk}^* + k_{ik}\Delta \psi_k \qquad (32)$$

where $\Delta V_{PCCDQk} = [\Delta V_{PCCDk}, \Delta V_{PCCQk}]^T$, and $A_{PCCk}$ is the parameter matrix.

*3) Network and load models within MG$_k$:* The network and load models [9] are developed based on the lumped, series RL feeder lines and the RL-type constant-impedance loads, respectively, i.e.,

$$\Delta \dot{i}_{lineDQk} = A_{netk}[\Delta i_{lineDQk}] + B_{netk}[\Delta v_{bDQk}] + C_{netk}\Delta \omega_g \qquad (33)$$

$$\Delta \dot{i}_{loadDQk} = A_{loadk}[\Delta i_{loadDQk}] + B_{loadk}[\Delta v_{bDQk}] + C_{loadk}\Delta \omega_g \qquad (34)$$

where $\Delta i_{lineDQk}$, $\Delta i_{loadDQk}$ and $\Delta v_{bDQk}$ are variables of all lines, loads and bus voltages within MG$_k$, respectively. The deviation of $i_{oki}$ of all the DG units and $i_{PCCk}$ of MG$_k$, Fig. 2, is denoted as $\Delta i_{oDQk}$ and $\Delta i_{PCCDQk}$ in the common DQ-frame, respectively. Then $\Delta v_{bDQk}$ is represented as [9]

$$\Delta v_{bDQk} = R_N(M_{DGk}[\Delta i_{oDQk}] + M_{netk}[\Delta i_{lineDQk}] + M_{loadk}[\Delta i_{loadDQk}] + M_{PCCk}[\Delta i_{PCCDQk}]) \qquad (35)$$

Since PCC$_k$ is also a bus within MG$_k$, $\Delta V_{PCCDQk}$ can be expressed in terms of $\Delta i_{oDQk}$, $\Delta i_{lineDQk}$, $\Delta i_{loadDQk}$, $\Delta i_{PCCDQk}$ and named as $\Delta V_{PCCDQk}$ *expression.*

*4) The complete model of MG$_k$:* Use the state variables of all DG units within MG$_k$ as $\Delta X_{DGk}$, then combine (28), (31), (33), (34) and replace $\Delta V_{fk}^*$, $\Delta v_{bDQk}$, $\Delta V_{PCCDQk}$ with (32), (35) and $\Delta V_{PCCDQk}$ expression, respectively. The small-signal dynamic model of MG$_k$ is obtained, i.e.,

$$[\Delta \dot{S}_{MGk}] = A_{MGk}[\Delta S_{MGk}] + B_{MGk}\Delta V_{PCCk}^* + C_{MGk}[\Delta i_{PCCDQk}] \qquad (36)$$

where $\Delta S_{MGk} = [\Delta X_{DGk}, \Delta i_{lineDQk}, \Delta i_{loadDQk}, \Delta \psi_k]^T$, $A_{MGk}$, $B_{MGk}$ and $C_{MGk}$ are parameter matrices.

### B. NMG Layer Modeling

The NMG layer modeling covers the dynamics of TC and DQC controllers as well as MV lines and loads. Note that each MG will be viewed as a black box and referred to as a MG unit.

1) *MG unit model:* By linearizing (19)-(22), (24)-(25) around an operating point, the model of $MG_k$ becomes

$$[\Delta \dot{X}_{\text{MG}k}] = A_{\text{MG}k}[\Delta X_{\text{MG}k}] + B_{\text{MG}k}[\Delta v_{\text{bDQ}k}] + C_{\text{MG}k}\Delta \omega_{\text{g}} +$$
$$\sum_{l \in \mathcal{H}_k} F_{\text{MG}k}[\Delta X_{\text{MG}l}] + H_{\text{MG}k}\Delta V_f^* + I_{\text{MG}k}[\Delta V_{\text{PCCDQ}k}] \quad (37)$$

where $\Delta v_{\text{bDQ}k}$ is the deviation of MV bus voltage $v_{\text{b}k}$ in the common DQ-frame, $A_{\text{MG}k}$, $B_{\text{MG}k}$ $C_{\text{MG}k}$, $F_{\text{MG}k}$ and $H_{\text{MG}k}$ are parameter matrices. Note that $F_{\text{MG}k}$ reflects the correlation between unit $MG_k$ and its neighbors $MG_l$, $l \in \mathcal{H}_k$. The state variables of each MG unit are

$$[\Delta X_{\text{MG}k}] =$$
$$\left[\Delta \delta_k, \Delta P_{\text{PCC}k}, \Delta Q_{\text{PCC}k}, \Delta \Omega_k, \Delta \lambda_k, \Delta h_k, \Delta i_{\text{PCC}dk}, \Delta i_{\text{PCC}qk}\right]^T \quad (38)$$

2) *Critical bus voltage controller model:* Denote $\psi$ as the state of (23), i.e.,

$$\dot{\psi} = V_c^* - V_c \quad (39)$$

where $V_c = \sqrt{V_{\text{cD}}^2 + V_{\text{cQ}}^2}$. By linearizing (23) and (39), the model of the critical bus voltage controller can be obtained,

$$\Delta \dot{\psi} = -A_c[\Delta V_{\text{cDQ}}] \quad (40)$$
$$\Delta V_f^* = -k_p A_c[\Delta V_{\text{cDQ}}] + k_i \Delta \psi \quad (41)$$

where $\Delta V_{\text{cDQ}} = [\Delta V_{\text{cD}}, \Delta V_{\text{cQ}}]^T$, and $A_c$ is the parameter matrix.

3) *MV network and load models:* The modeling of MV network and load is the same as that in MG layer modeling and can be expressed as

$$\Delta \dot{i}_{\text{lineDQ}} = A_{\text{net}}[\Delta i_{\text{lineDQ}}] + B_{\text{net}}[\Delta v_{\text{bDQ}}]$$
$$+ C_{\text{net}}\Delta \omega_{\text{g}} \quad (42)$$

$$\Delta \dot{i}_{\text{loadDQ}} = A_{\text{load}}[\Delta i_{\text{loadDQ}}] + B_{\text{load}}[\Delta v_{\text{bDQ}}]$$
$$+ C_{\text{load}}\Delta \omega_{\text{g}} \quad (43)$$

where $\Delta i_{\text{lineDQ}}$, $\Delta i_{\text{loadDQ}}$ and $\Delta v_{\text{bDQ}}$ are variables of MV lines, loads and buses, respectively. $\Delta i_{\text{PCCDQ}}$ denotes $\Delta i_{\text{PCCDQ}k}$ of all the MG units. Then $\Delta v_{\text{DQ}}$ is represented as

$$\Delta v_{\text{bDQ}} = R_N(M_{\text{MG}}[\Delta i_{\text{PCCDQ}}] + M_{\text{net}}[\Delta i_{\text{lineDQ}}]$$
$$+ M_{\text{load}}[\Delta i_{\text{loadDQ}}]) \quad (44)$$

Since the critical bus is also a MV bus, $\Delta v_{\text{cDQ}}$ can be expressed in terms of $\Delta i_{\text{PCCDQ}}$, $\Delta i_{\text{lineDQ}}$ and $\Delta i_{\text{loadDQ}}$ and named as $\Delta V_{\text{cDQ}}$ *expression.*

4) *Complete the NMG layer model:* Denote the state variables of all MG units in the NMG system as $\Delta X_{MG}$. Combine (37), (40), (42), (43) and replace $\Delta V_f^*$, $\Delta v_{bDQ}$, $\Delta V_{cDQ}$ with (41), (44) and $\Delta V_{cDQ}$ *expression.* Then, the small-signal dynamic model of the NMG layer is obtained. That is,

$$[\Delta \dot{S}_{\text{NMG}}] = A_{\text{NMG}}[\Delta S_{\text{NMG}}] + B_{\text{NMG}}\Delta V_{\text{PCCDQ}k} \quad (45)$$

where $\Delta S_{\text{NMG}} = [\Delta X_{\text{MG}}, \Delta i_{\text{lineDQ}}, \Delta i_{\text{loadDQ}}, \Delta \psi]^T$, $A_{\text{NMG}}$, $B_{\text{NMG}}$ are parameter matrices.

### C. Complete NMG System Model

In (36) and (45), the coupling states $\Delta V_{\text{PCC}k}^*$, $\Delta i_{\text{PCCDQ}k}$ and $\Delta V_{\text{PCCDQ}k}$ can be dealt with as follows: i) linearize (24), then $\Delta V_{\text{PCC}k}^*$ in (36) can be represented by $\Delta X_{\text{MG}k}$ in (38) which is part of $\Delta S_{\text{NMG}}$; ii) represent $\Delta i_{\text{PCCDQ}k}$ by $\Delta X_{\text{MG}k}$ which is part of $\Delta S_{\text{NMG}}$ and iii) replace $\Delta V_{\text{PCCDQ}k}$ by $\Delta V_{\text{PCCDQ}k}$ *expression* then represented by $\Delta S_{\text{MG}}$. By combining $m$ MG layer models (36) and NMG layer model (45), the complete NMG system model can be obtained as

$$[\Delta \dot{S}_{\text{sys}}] = A_{\text{sys}}[\Delta S_{\text{sys}}] \quad (46)$$

where $\Delta S_{\text{sys}} = [\Delta S_{\text{MG1}}, \dots, \Delta S_{\text{MG}m}, \Delta S_{\text{NMG}}]^T$.

## V. NUMERICAL STUDY

To validate the effectiveness of the proposed two-layer distributed control method, stability analyses and time-domain simulation studies in the PSCAD/EMTDC platform are carried out in this section based on a test NMG system.

### A. Test System

The test NMG system including 3 MGs is shown in Fig. 4. The circuit breaker (CB) 1, 2 and 3 are closed. The rated voltages of MV and LV network are 10kV and 0.38kV, respectively. Each MG connects with the MV feeder through a 10kV/0.38kV $\Delta$/$Y_g$ transformer. $L_c$ is the coupling inductance. Tables I, II and III provide the system and controller parameters, respectively.

TABLE I. ELECTRICAL PARAMETERS OF THE NMG SYSTEM

| Line | Line1: $0.08 + \text{j}0.12$ $\Omega$ , Line2: $0.05 + \text{j}0.07$ $\Omega$ , Line3: $0.07 + \text{j}0.11$ $\Omega$ , Line4,8,10,11: $0.15 + \text{j}0.05$ $\Omega$ , Line6,7,12: $0.11 + \text{j}0.07$ $\Omega$ , Line5,9: $0.11 + \text{j}0.11$ $\Omega$ , |
|---|---|
| Load | Load9=100kw+30kvar, Load10=20kw+5kvar, Load1,5,8=15kw+7.5kvar, Load3,6=12kw+5kvar, Load2=40kw+15kvar, Load4,7=50kw+20kvar |
| Transformer | T1/T2/T3: 1Mva, $u_k = 4\%$, $r_k = 1\%$, 10/0.38kV($\Delta$/$Y_g$) |

TABLE II. PARAMETERS OF PCS AND TCS

| Parameters | DG11 DG12 DG13 DG31 DG32 DG33 | DG21 DG22 DG23 | MG1 MG3 | MG2 |
|---|---|---|---|---|
| $D_{Pk\ell}/D_{Pk}$ (Hz/kW $\cdot 10^{-3}$) | 16.67 | 20 | 8.33 | 10 |
| $D_{Qk\ell}/D_{Qk}$ (kV/kvar $\cdot 10^{-3}$) | 0.78 | 0.52 | 0.39 | 0.26 |
| $P_{max}/P_{sMG}$ (kW) | 60 | 50 | 120 | 100 |
| $Q_{max}/Q_{sMG}$ (kvar) | 20 | 30 | 40 | 60 |

TABLE III. PARAMETERS OF DSCS AND DQCS

| Parameters | DSC level | DQC level |
|---|---|---|
| $c_{\omega k\ell}/c_{\omega k}$ | 560 | 2500 |
| $c_{pk\ell}/c_{pk}$ | 50 | 2500 |
| $c_{vk\ell}/c_{vk}$ | 20 | 12 |
| $c_{qk\ell}/c_{qk}$ | 100 | 25 |
| $k_{pk}/k_p$ | 30 | 5 |
| $k_{ik}/k_i$ | 0.05 | 20 |

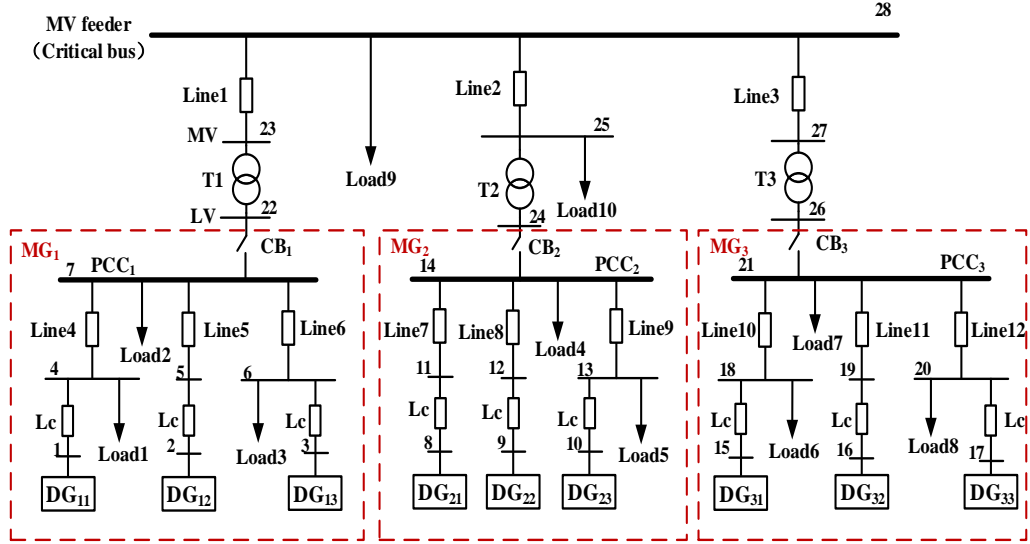

Fig. 4  A schematic diagram of the test NMG system

Reference values of the system frequency and critical bus voltage are given as $f_{sys}^* = 50$Hz and $V_c^* = 1$p.u., respectively. The communication networks $G_k$ and $\tilde{G}$ in DSC and DQC level are assumed to have the same topology as shown in Fig. 5. Moreover, in the communication network $G_k$ of MG$_k$, only the node DG$_{k1}$ receives reference value $(\omega_{MGk}/V_{tk}^*)$ with pinning gain $g_{k1} = 1$. Similarly, in $\tilde{G}$, only the MG$_1$ receives reference value $(\omega_{sys}^*/V_t^*)$ with pinning gain $g_1 = 1$.

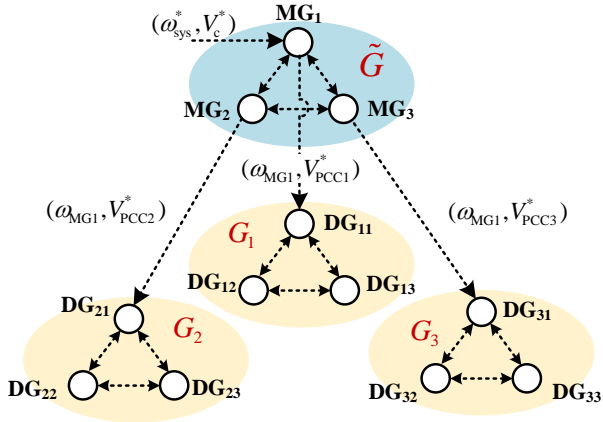

Fig. 5  Topology of two-layer distributed communication network

### B. Eigen-Analysis Results

The system dynamics and stability analysis results with the proposed control method are presented in this subsection.

*1) Participation factors:* Fig. 6 compares the low-frequency eigenvalue spectrum of single MG$_1$, MGs$_2$, MG$_3$ and the NMG system. Note that a single MG only employs the MG layer controllers, in Section III, with $\omega_{MGk}$ and $V_{PCCk}^*$ set as desired values. The typical dominant modes of the NMG system are labelled as mode $i$ ($i = 1,2,...,8$).

Fig. 6 indicates that interconnecting MGs (i) significantly changes the shaping of the eigenvalues on the complex plane, (ii) still retains the single MG modes, i.e. modes 1-3 and (iii) introduces four pairs of low-damping modes (modes 4-7) leading to much more oscillatory system responses as compared with single MGs.

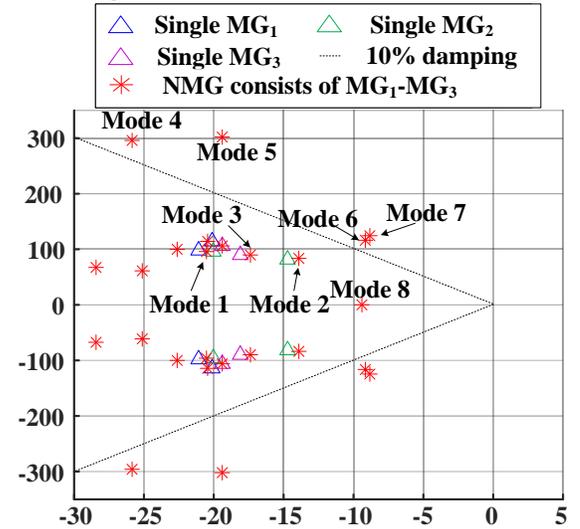

Fig. 6  Low-frequency eigenvalue spectrum of three single MGs and NMG system

To identify the correlation between system states and dominant oscillatory modes, a participation factor analysis is conducted. Fig. 7 illustrates the participation factors of MG layer states ($\Delta\delta_{ki}, \Delta P_{ki}, \Delta Q_{ki}, \Delta\Omega_{ki}, \Delta\lambda_{ki}, \Delta h_{ki}$), NMG layer states ($\Delta\delta_k, \Delta P_k, \Delta Q_k, \Delta\Omega_k, \Delta\lambda_k, \Delta h_k$) and the critical bus voltage controller state ($\Delta\psi$). As indicated by Fig. 7, mode 1 is a typical MG inner mode which is almost only affected by states of DG$_{11}$~DG$_{13}$ units within MG$_1$. Mode 5 and 7 are inter-coupling modes mainly affected by states of both MG and NMG layers. For simplicity, the participation factors of modes 2, 3 (MG inner modes of MG$_2$ and MG$_3$) and modes 4, 6 (inter-coupling modes) are not presented. Besides, mode 8 is mainly affected by the critical bus voltage controller. The strongly associated states, controllers and parameters with modes 1-8 are summarized in Table IV.

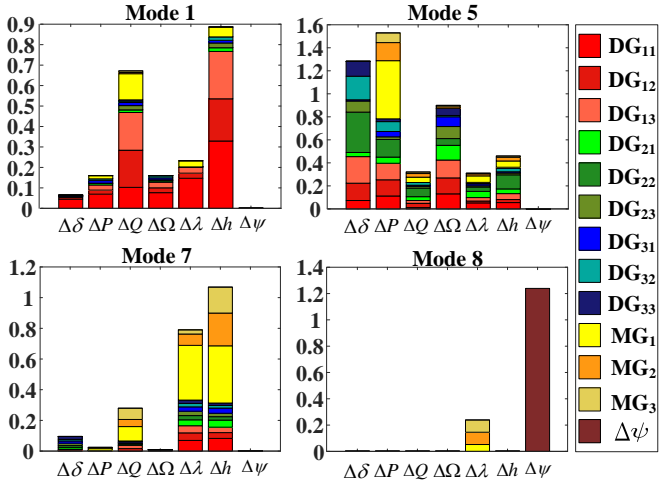

Fig. 7  Participation factors of modes 1, 5, 7 and 8

Table IV indicates that the most dominant modes 4-7 (with damping less than 10%) are affected by control parameters of DSCs and DQCs. Therefore, the impact of these parameters on system stability should be carefully analyzed.

*2) Sensitivity analysis of DSC and DQC parameters*: Fig. 8 shows the traces of modes 4-7 as a function of $c_{pki}$ and $c_{qk}$. The impacts of other DSC and DQC parameters are summarized in Table IV. Fig. 8 shows that the variation of parameters may bring instability risk to the system. The summary in Table IV indicates that a dominant mode may be affected by multiple controllers and their parameters. Besides, a control parameter may affect different modes.

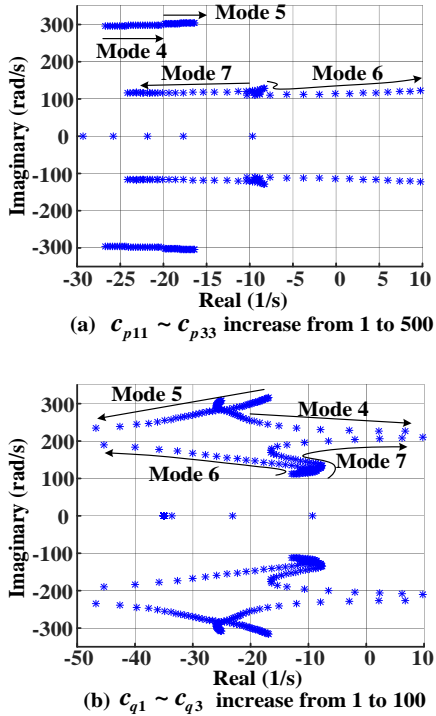

Fig. 8  Traces of the most dominant modes 4-7.

*3) Summary*: The above results reveal that (i) interconnecting MGs introduces new low-frequency oscillatory modes and therefore complicates the system dynamic behavior; (ii) the new low-damping modes (modes 4-7) reduce the stability margin due to the coupling among neighboring MGs and between the two control layers; (iii) the dominant oscillatory modes are affected by multiple controllers. Note that the parameters in Table III are carefully tuned based on the guidelines in Table IV.

### C. Time-Domain Simulation Results in PSCAD/EMTDC

Three cases are designed for the simulation. Case 1 demonstrates the steady-state performance of the proposed control strategy, namely the capability to meet control objectives (i)~(iv) simultaneously under normal condition. Case 2 verifies the system dynamic performance under communication failures as well as sudden load changes. Case 3 verifies the plug-and-play functionality of DG and MG units.

*1) Case 1-steady-state performance:* The PCs are initially engaged and the DSCs and TCs are activated at t=1.5s and DQCs are employed at t=3s. Fig. 9 (a) indicates that a 0.25 Hz frequency deviation is introduced by TCs while the DQCs restore the system frequency to 50Hz (objective i). Fig. 9 (b) indicates that the critical bus voltage is restored to 1 p.u. by the DQC (objective ii). Fig. 9 (d) indicates that after t=3s the DQCs achieve accurate reactive power sharing among MGs with the ratios of $Q_1$ to $Q_3$ being 2:3:2 (objective iii). Fig. 9 (f) indicates that after t=1.5s the DSCs realize accurate reactive power sharing among DGs within each MG with the ratios of $Q_{k1}$ to $Q_{k3}$ being 1:1:1 (objective iv). Fig. 9 (c) and (e) indicate that the output active powers of MGs and DGs are always accurately shared. Note that the power sharing ratios are presented in Table II.

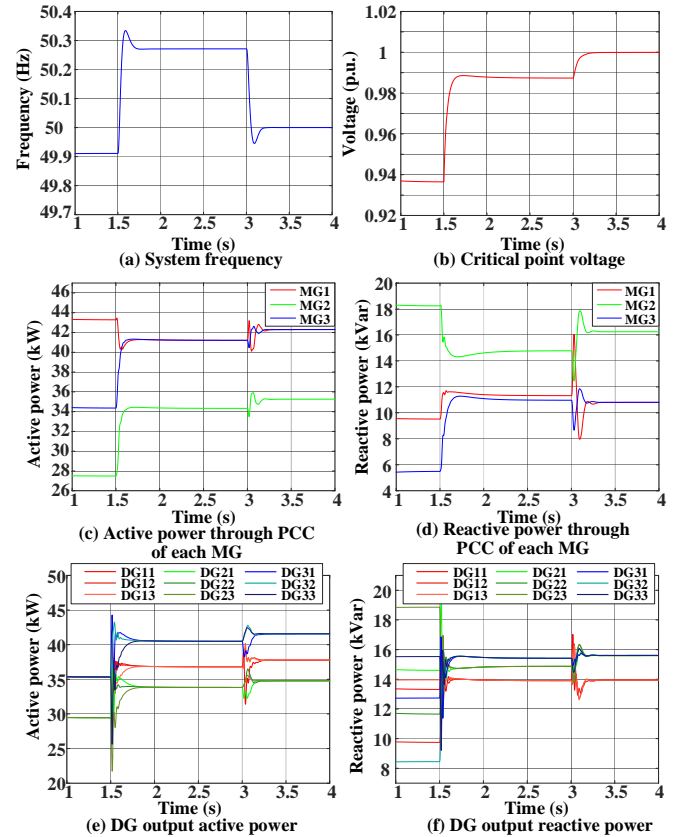

Fig. 9  Steady-state performance of the proposed method.

TABLE IV. RESULTS OF PARTICIPATION FACTORS ANALYSIS

| Modes | Strongly associated control layer | Strongly associated states | Associated controllers | Parameters of the associated controller | Impact of parameters on mode damping (MD) |
|-------|-----------------------------------|----------------------------|------------------------|------------------------------------------|-------------------------------------------|
| 1-3 | MG layer-DSClevel | $\Delta Q_{ki}, \Delta \lambda_{ki}, \Delta h_{ki}$ | (15),(25) | $c_{vki}, c_{qki}$ | $c_{vki}, c_{qki} \uparrow$ (MD $\downarrow$) |
| 4-5 | MG layer-DSClevel | $\Delta Q_{ki}, \Delta \lambda_{ki}, \Delta h_{ki}$ | (11),(15) | $c_{vki}, c_{qki}$ | Mode 4: $c_{vki}, c_{qki}, c_{vk}, c_{qk} \uparrow$ (MD $\downarrow$); |
| | NMG layer-DQClevel | $\Delta Q_k, \Delta \lambda_k, \Delta h_k$ | (22),(25) | $c_{vk}, c_{qk}$ | Mode 5: $c_{vki}, c_{vk} \uparrow$ (MD $\downarrow$);  $c_{qki}, c_{qk} \uparrow$ ( MD $\uparrow$) |
| 6-7 | MG layer-DSClevel | $\Delta \delta_{ki}, \Delta P_{ki}, \Delta \Omega_{ki}$ | (9) | $c_{\omega ki}, c_{pki}$ | Mode 6: $c_{pki} \uparrow$ (MD $\downarrow$);   $c_{\omega ki}, c_{qk} \uparrow$ (MD $\uparrow$) |
| | NMG layer-DQClevel | $\Delta h_k$ | (25) | $c_{qk}$ | Mode 7: $c_{qk} \uparrow$ (MD $\downarrow$);   $c_{\omega ki}, c_{pki} \uparrow$ (MD $\uparrow$) |
| 8 | NMG layer-DQClevel | $\Delta \psi$ | (23) | $k_p, k_i$ | $k_p, k_i \uparrow$ (MD $\downarrow$);   $k_p, k_i \downarrow$ (MD $\uparrow$) |

*2) Case 2- communication link failures:* In this case, all the controllers are activated at t=0.8s and then the system reaches a steady state. It is worth to be noted that, based on the proof in [34], for tracking synchronization problem and regulator synchronization problem, if a spanning tree exists in the corresponding distributed communication network and $g_k \neq 0$ for at least one root node, the proposed controllers can reach a steady state and objectives (i)~(iv) can still be realized.

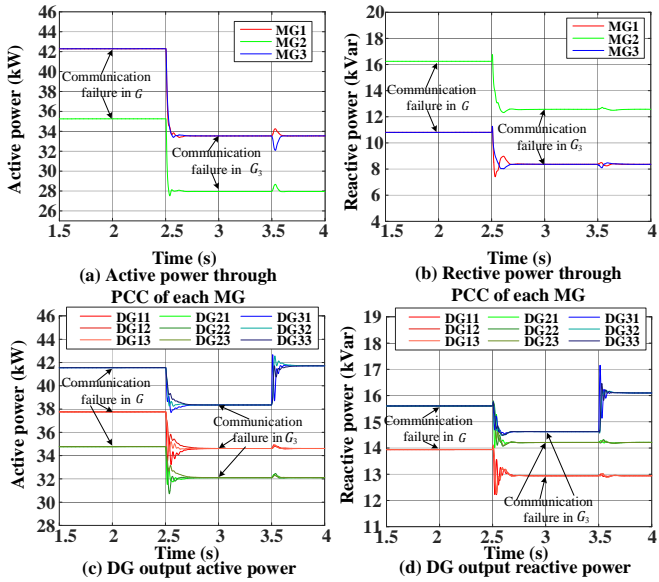

(a) Active power through PCC of each MG

(b) Rective power through PCC of each MG

(c) DG output active power

(d) DG output reactive power

Fig. 10 System behaviors when communication failures happen in $\tilde{G}$ and $G_3$.

Stage 1 (1.5s-3s): in the upper communication network $\tilde{G}$, the communication link between $MG_2$ and $MG_3$, as shown in Fig. 5, fails at $t = 2$ s. Subsequently, 25% of load 9 is switched off at t=2.5s. The results in Fig. 10 show that the steady-state objectives (i)-(iv) can still be achieved after the communication link failure since the remaining communication network still contains a spanning tree.

Stage 2 (3s-4s): during this stage, a worse scenario which refers to communication failures happening in both upper network $\tilde{G}$ and lower network $G_k$ is set up. After one communication link fails at t=2s in upper communication network $\tilde{G}$, for lower communication network $G_3$ in $MG_3$, the communication link between $DG_{32}$ and $DG_{33}$, as shown in Fig. 10, fails at $t = 3$s. Subsequently, at $t = 3.5$s, 80% of load 6 which is the internal load of $MG_3$ is switched on. The results in Fig. 10 shows that the steady-state objectives (i)-(iv) can be realized since both the network $\tilde{G}$ and $G_3$ still contain a spanning tree after communication failure in this stage. Besides, Fig. 10 also indicates that the NMG system reaches a steady state within 0.5s after a disturbance and no significant overshoot is observed even under communication failure events.

*3) Case 3- plug-and-play operation:* In this study, all the controllers are activated at t=0.8s and then the system reaches a steady state.

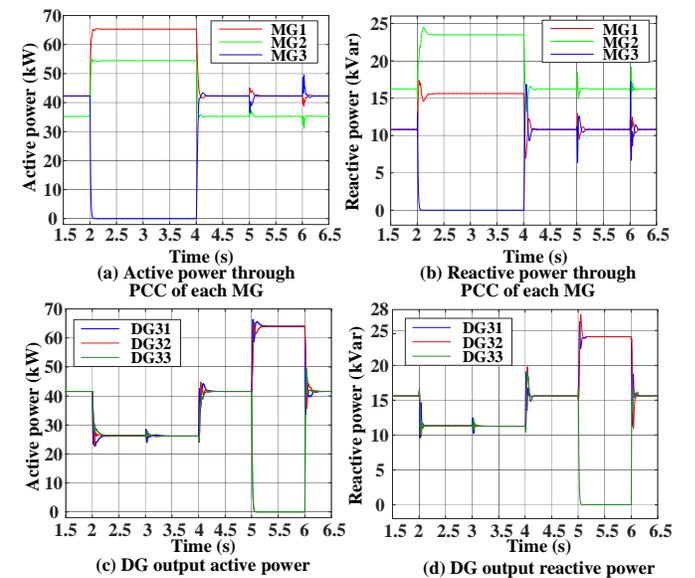

(a) Active power through PCC of each MG

(b) Reactive power through PCC of each MG

(c) DG output active power

(d) DG output reactive power

Fig. 11 System behaviors under plug-and-play operation of $MG_3$ and $DG_{33}$.

Stage 1 (2s-4s): $MG_3$ is disconnected at t=2s and reconnected at t=4s to evaluate the plug and play capability of $MG_3$. Note that (i) the MG will lose all the communication links with its neighbouring units when it disconnects with the NMG system, then these links will recover after its reconnection; (ii) the synchronization process is necessary for $MG_3$ before its reconnection (specifically, in this study, the synchronization of $MG_3$ starts at t=3s during its islanded state); (iii) the MGs will transfer to the islanded operation state with only MG layer controllers employed after the disconnection event at t=2s, and the reference angular frequency $\omega_{MGk}$ and reference PCC voltage $V^*_{PCCk}$ will be set as the rated value $2 * pi * 50$ rad/s and 1.0 p.u, respectively, to maintain a stable operation of the islanded $MG_3$. The active and reactive power through PCC of each MG are presented in Fig. 11 (a) and (b), respectively.

Stage 2 (5s-6s): DG$_{33}$ in MG$_3$ disconnects at t=5s and reconnects at t=6s to evaluate the plug and play capability of DGs. Similarly, DG$_{33}$ will lose all the communication links with its neighboring units when it disconnects and the communication links will recover after its reconnection. The synchronization process of DG$_{33}$ starts immediately after it disconnects at t=5s. The active and reactive power of DG units in MG$_3$ are presented in Fig. 11 (c) and (d), respectively.

The above results show that after the disconnection of MG$_3$ in the NMG-control layer and DG$_{33}$ in the MG-control layer, the objectives (i)~(iv) can still be realized (for the sake of simplicity, only the results of active and reactive power are given here). This is because of the remaining communication network $\tilde{G}$ of NMG-control layer and $G_3$ of MG-control layer still contain a spanning tree. After MG$_3$ and DG$_{33}$ reconnect, their power sharing objectives can be realized. Besides, during plug-and-play operation, no significant overshoots are observed, and the time of recovering to a steady state is within 0.4s.

## VI. EXPERIMENTAL VALIDATION

This section provides experimental results to validate the practical implementation feasibility of the proposed methods.

### A. Experimental setup

Fig.12 shows the experimental NMG system setup, which includes a real-time dSPACE 1006 platform, four Danfoss inverters, resistive loads, inductive loads, line impedances, switches and the control desk. The control strategy is programmed and executed in the dSPACE 1006 platform to switch the inverters. The switching frequency of the inverters is set to be 10kHz.

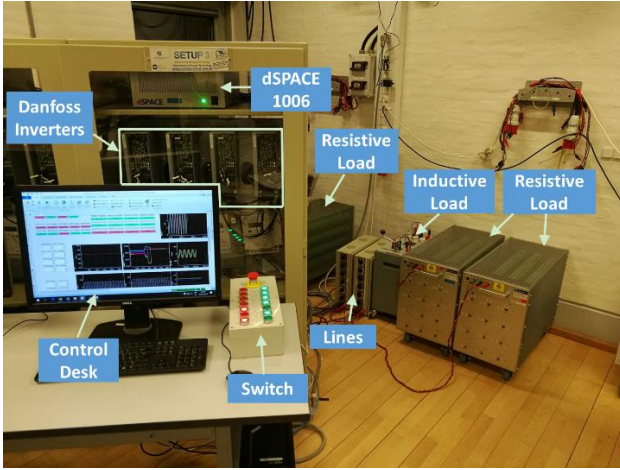

Fig. 12  Experimental setup in the laboratory

The physical configuration of the experimental NMG system is shown in Fig.13. There are two microgrids in the system and each microgrid consists of two DG units. The circuit breakers CB$_1$ and CB$_2$ are closed. MG$_1$ and MG$_2$ are connected to the critical bus through line impedances. The system rated frequency is 50Hz and the rated rms voltage is 200V.

Table V provides the electrical parameters of the experimental NMG system. Parameters of the four level controllers are shown in Table VI and VII. From Table VI, it can be seen that (i) the active and reactive power capacities between MG$_1$ and MG$_2$ are equal, i.e., $P_{sMG1} = P_{sMG2}$ and $Q_{sMG1} = Q_{sMG2}$, and (ii) the ratios of the active and reactive power capacities between DG$_{k1}$ and DG$_{k2}$ are both 4:3, i.e., $P_{maxk1}:P_{maxk2} = 4:3$ and $Q_{maxk1}:Q_{maxk2} = 4:3$, $k = 1, 2$.

Fig.14 shows the topology of two-layer distributed communication networks for the experimental NMG system. In the upper network $\tilde{G}$, MG$_1$ receives reference values $\omega_{sys}^*$ and $V_c^*$. In the lower network $G_1$ and $G_2$, DG$_{11}$ and DG$_{21}$ receive the references from $\tilde{G}$.

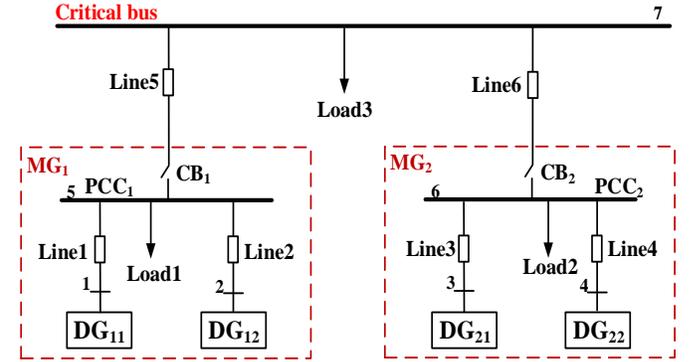

Fig. 13  Physical configuration of experimental NMG system

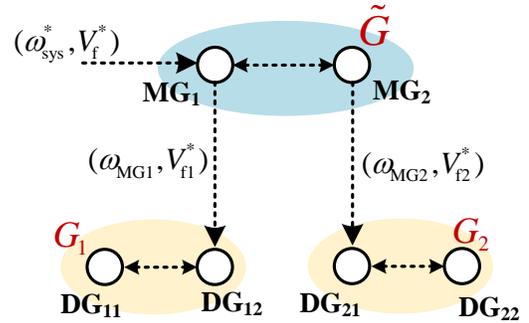

Fig. 14  Topology of the two-layer distributed communication network for experimental NMG system

TABLE V.  ELECTRICAL PARAMETERS OF EXPERIMENTAL NMG SYSTEM

| Inverter | Filter inductance: 1.8mH<br>Filter capacitance: 27µF |
|---|---|
| Line | Line1=1.8mH, Line2=1.8mH,<br>Line3=1.8mH, Line4=1.8mH,<br>Line5=1.9Ω+2.5mH, Line6=1.6Ω+2.1mH |
| Load | Load1=92Ω, Load2=153.3Ω,<br>Load3=38.1+j32.9Ω |

TABLE VI.  PARAMETERS OF PCS AND TCS FOR EXPERIMENTAL NMG SYSTEM

| Parameters | DG$_{11}$ | DG$_{12}$ | DG$_{21}$ | DG$_{22}$ | MG$_1$ | MG$_2$ |
|---|---|---|---|---|---|---|
| $D_{Pki}/D_{Pk}$ <br> (Hz/W · $10^{-3}$) | 0.625 | 0.833 | 0.625 | 0.833 | 0.357 | 0.357 |
| $D_{Qki}/D_{Qk}$ <br> (V/Var · $10^{-3}$) | 6.479 | 8.639 | 6.479 | 8.639 | 3.702 | 3.702 |
| $P_{max}/P_{sMG}$ <br> (kW) | 1.8 | 1.35 | 1.8 | 1.35 | 3.15 | 3.15 |
| $Q_{max}/Q_{sMG}$ <br> (kvar) | 1.2 | 0.9 | 1.2 | 0.9 | 2.1 | 2.1 |

| Parameters | DSC level | DQC level |
|---|---|---|
| $c_{\omega ki}/c_{\omega k}$ | 400 | 80 |
| $c_{pki}/c_{pk}$ | 400 | 80 |
| $c_{vki}/c_{vk}$ | 150 | 30 |
| $c_{qki}/c_{qk}$ | 20 | 2 |
| $k_{pk}/k_p$ | 1.2 | 0.3 |
| $k_{ik}/k_i$ | 42 | 10 |

### B. Experimental results

The system is initially operated with PC. At t=6.9s the DSC and TC are activated, and at t=8.55s the DQC is activated. Fig.15 shows the corresponding experimental results. Fig.15(a) and (b) indicate that after t=8.55s the DQC can restore the system frequency and critical bus voltage to their rated values 50Hz and 1 p.u., i.e., objectives (i) and (ii) are realized. Fig.15(c) indicates that the active powers through PCC of each MG are equal after applying TC t=6.9s (objective (iii)-(1)). Fig.15(d) indicates that the reactive powers through PCC of each MG are equal after applying DQC at t=8.55s (objective (iii)-(2)). Fig.15(e) and (f) indicate that after applying DSC at t=6.9s, the output active power and reactive power of DG units in each MG can realize accurate sharing with $P_{k1}$ : $P_{k2}$ and $Q_{k1}$ : $Q_{k2}$ both being 4:3 (objective (iv)).

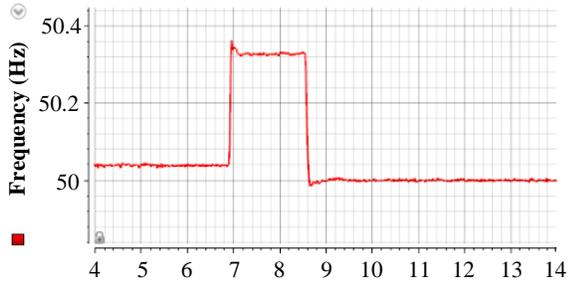

(a) system frequency

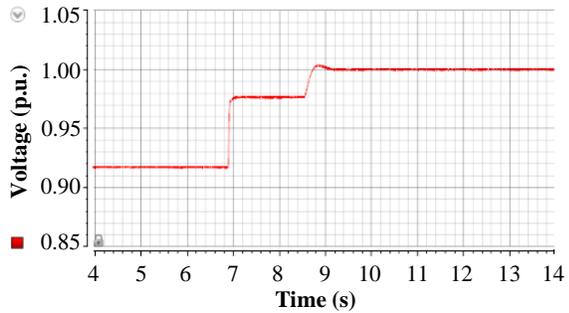

(b) critical bus voltage

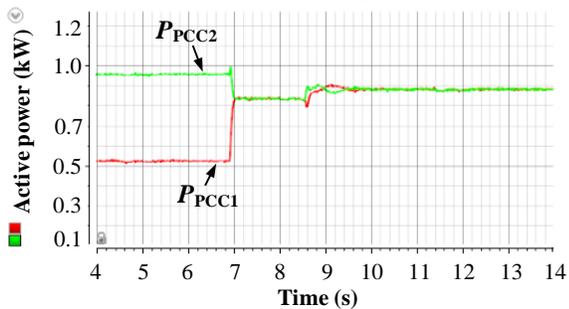

(c) Active power through PCC of each MG

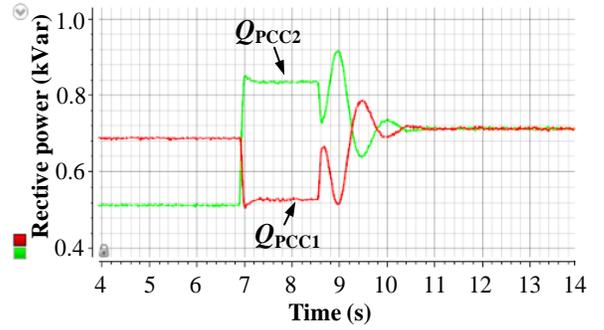

(d) Reactive power through PCC of each MG

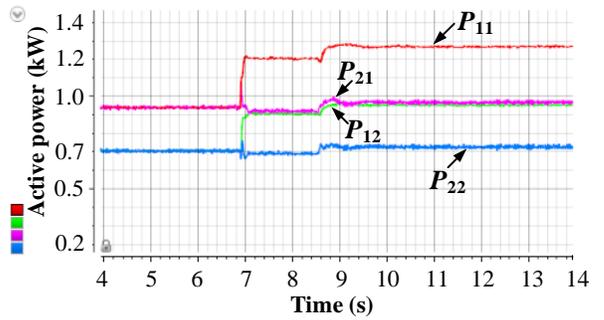

(e) DG output active power

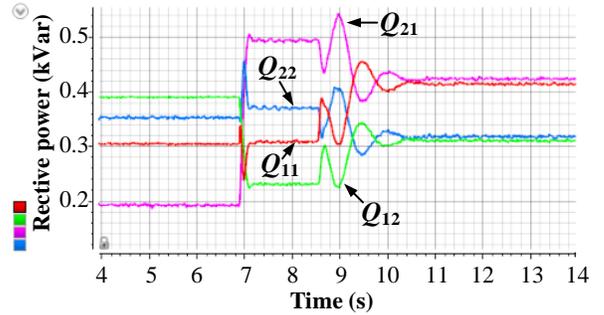

(f) DG output reactive power

Fig. 15  Experimental results (DSC and TC are activated at t=6.9s, and DQC is activated at t=8.55s)

## VII. CONCLUSION

This paper presents a two-layer, four-level distributed control method for islanded NMG systems. Under the proposed control architecture, the proprietary information of each MG is well protected. The proposed method is capable of simultaneously meeting objectives including the regulation of frequency and critical bus voltage, as well as accurate power sharing in both the MG and NMG layers. Moreover, small-signal analysis with the proposed control method indicates that i) interconnecting MGs results in low-damping inter-coupling modes which may lead to system instability, ii) the system dominant modes are strongly affected by controller parameters, especially those in the DSC and DQC levels. Time-domain simulations as well as experimental results on a NMG test system validate the effectiveness of the proposed methods.

## IX. Biographies

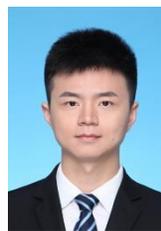
**Xiaoyu Wu** (S'15) received his B.S degree in electrical engineering from Beijing Jiaotong University, Beijing, China, in 2013.

He is currently pursuing his Ph.D. degree at Beijing Jiaotong University, Beijing, China. During 2015-2016, he was a visiting scholar at the North Carolina State University, Raleigh, NC, USA. His research interests include demand response, control and stability analysis of microgrid.

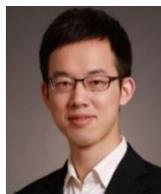
**Yin Xu** (S'12–M'14-SM'18) received his B.E. and Ph.D. degrees in electrical engineering from Tsinghua University, Beijing, China, in 2008 and 2013, respectively.

He is currently a Professor with the School of Electrical Engineering at Beijing Jiaotong University, Beijing, China. During 2013–2016, he was an Assistant Research Professor with the School of Electrical Engineering and Computer Science at Washington State University, Pullman, WA, USA. His research interests include power system resilience, distribution system restoration, power system electromagnetic transient simulation, and AC-DC hybrid power systems.

Dr. Xu is currently serving as Secretary of the Distribution Test Feeder Working Group under the IEEE PES Distribution System Analysis Subcommittee.

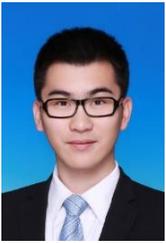

**Xiangyu Wu** (S'13–M'17) received the B.S. degree from the Department of Electrical Engineering, Zhejiang University, Hangzhou, China, in 2012 and the Ph.D. degree in electrical engineering from Tsinghua University, Beijing, China, in 2017.

He is currently a postdoctoral researcher at the School of Electrical Engineering, Beijing Jiaotong University, Beijing, China. In 2015, he was a visiting scholar at University of Toronto, Toronto, ON, Canada. Since 2018, he is a guest researcher at the Department of Energy Technology, Aalborg University, Aalborg, Denmark. His research interests include the control and stability analysis of microgrid and power electronic based power systems.

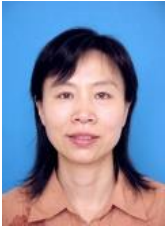

**Jinghan He** (M'07-SM'18) received the M.Sc. degree in electrical engineering from the Tianjin University, Tianjin, China, in 1994 and the Ph.D. degree in electrical engineering from Beijing Jiaotong University, Beijing, China, in 2007.

She is currently a Professor at Beijing Jiaotong University, Beijing, China. Her main research interests are protective relaying, fault distance measurement, and location in power systems.

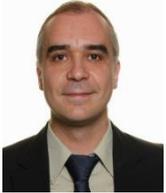

**Josep M. Guerrero** (S'01-M'04-SM'08-FM'15) received the B.S. degree in telecommunications engineering, the M.S. degree in electronics engineering, and the Ph.D. degree in power electronics from the Technical University of Catalonia, Barcelona, in 1997, 2000 and 2003, respectively. Since 2011, he has been a Full Professor with the Department of Energy Technology, Aalborg University, Denmark, where he is responsible for the Microgrid Research Program (www.microgrids.et.aau.dk). From 2014 he is chair Professor in Shandong University; from 2015 he is a distinguished guest Professor in Hunan University; and from 2016 he is a visiting professor fellow at Aston University, UK, and a guest Professor at the Nanjing University of Posts and Telecommunications.

His research interests is oriented to different microgrid aspects, including power electronics, distributed energy-storage systems, hierarchical and cooperative control, energy management systems, smart metering and the internet of things for AC/DC microgrid clusters and islanded minigrids; recently specially focused on maritime microgrids for electrical ships, vessels, ferries and seaports. Prof. Guerrero is an Associate Editor for a number of IEEE TRANSACTIONS. He has published more than 450 journal papers in the fields of microgrids and renewable energy systems, which are cited more than 30,000 times. He received the best paper award of the IEEE Transactions on Energy Conversion for the period 2014-2015, and the best paper prize of IEEE-PES in 2015. As well, he received the best paper award of the Journal of Power Electronics in 2016. During five consecutive years, from 2014 to 2018, he was awarded by Thomson Reuters as Highly Cited Researcher. In 2015 he was elevated as IEEE Fellow for his contributions on "distributed power systems and microgrids."

**Chen-Ching Liu** (S'80–M'83–SM'90–F'94) received the Ph.D. degree from the University of California, Berkeley, CA, USA, in 1983.

He served as a Professor with the University of Washington, Seattle, WA, USA, from 1983–2005. During 2006–2008, he was Palmer Chair Professor with Iowa State University, Ames, IA, USA. During 2008–2011, he was a Professor and Acting/Deputy Principal of the College of Engineering, Mathematical and Physical Sciences with University College Dublin, Ireland. He was Boeing Distinguished Professor of Electrical Engineering and Director of the Energy Systems Innovation Center, Washington State University, Pullman, WA, USA. Currently, he is American Electric Power Professor and Director of the Center for Power and Energy, Virginia Polytechnic Institute and State University, Blacksburg, VA, USA. He is also a Research Professor at Washington State University.

Prof. Liu was the recipient of the IEEE PES Outstanding Power Engineering Educator Award in 2004. In 2013, he received the Doctor Honoris Causa from Polytechnic University of Bucharest, Romania. He served as Chair of the IEEE PES Technical Committee on Power System Analysis, Computing, and Economics during 2005–2006. He is a Fellow of the IEEE.

**Kevin P. Schneider** (S'00–M'06–SM'08) received his B.S. degree in Physics and M.S. and Ph.D. degrees in Electrical Engineering from the University of Washington, Seattle, WA, USA.

He is currently a Principal Research Engineer with the Pacific Northwest National Laboratory, working at the Battelle Seattle Research Center, Seattle, WA, USA. He is an Adjunct Faculty member with Washington State University, Pullman, WA, USA, and an Affiliate Assistant Professor with the University of Washington, Seattle, WA, USA. His main areas of research are distribution system analysis and power system operations.

Dr. Schneider is a licensed Professional Engineer in Washington State. He is the past Chair of the Distribution System Analysis Sub-Committee and current Secretary of the Power System Analysis, Computing, and Economics (PSACE) Committee.

**Dan T. Ton** received his B.S. degree is in electrical engineering and M.S. degree in business management from the University of Maryland, Baltimore, MD, USA.

He is a Program Manager of Smart Grid R&D within the U.S. Department of Energy (DOE) Office of Electricity Delivery and Energy Reliability (OE). He is responsible for developing and implementing a multi-year R&D program plan for next-generation smart grid technologies to transform the electric grid in the United States through public/private partnerships. Previously, he managed the Renewable Systems Integration program within the DOE Solar Energy Technologies Program.